\documentclass[11pt,a4paper]{amsart}

\usepackage{amsmath,amssymb,amsthm,mathtools,mathrsfs}
\usepackage[T1]{fontenc}
\usepackage{lmodern}
\usepackage{microtype}
\usepackage[a4paper,margin=2.45cm,headheight=15pt]{geometry}
\usepackage{booktabs,array}
\usepackage{enumitem}
\usepackage{xcolor}
\usepackage{hyperref}
\usepackage[nameinlink,capitalise,noabbrev]{cleveref}

\definecolor{paperblue}{HTML}{173A5E}
\definecolor{papergold}{HTML}{8A5A16}
\hypersetup{
	colorlinks=true,
	linkcolor=black,
	citecolor=papergold,
	urlcolor=paperblue,
	pdftitle={Genus one Coinvariants for Quasi-lisse  VOA},
	pdfauthor={Peng Shan, Wenbin Yan, Yongchang Zhu}
}

\allowdisplaybreaks[3]
\numberwithin{equation}{section}

 \newtheorem{theorem}{Theorem}[section]
 \newtheorem{prop}[theorem]{Proposition}
 \newtheorem{cor}[theorem]{Corollary}
 \newtheorem{lemma}[theorem]{Lemma}
 \newtheorem{defn}[theorem]{Definition}

 \newcommand{\la}{\langle}
 \newcommand{\ra}{\rangle}
\newcommand{\deff}{\stackrel{\mathrm{def}}{=}}

\DeclareMathOperator*{\Res}{Res}

\newcommand{\C}{\mathbb C}
\newcommand{\Z}{\mathbb Z}

\title[Elliptic Coinvariants of Quasi-lisse  Vertex Operator Algebras]{
 Elliptic Coinvariants  of Quasi-lisse  Vertex Operator Algebras}

\author{Peng Shan}
\address{P. Shan,  Department of Mathematical Sciences and Yau Mathematical Sciences Center, Tsinghua University, Beijing, 100084, China}
\email{pengshan@tsinghua.edu.cn}

\author{Wenbin Yan}
\address{W. Yan,  Yau Mathematical Sciences Center, Tsinghua University, Beijing, 100084, China}
\email{wbyan@tsinghua.edu.cn}

\author{Yongchang Zhu}
\address{Y. Zhu: Yau Mathematical Sciences Center and Beijing Institute of Mathematical Sciences and Applications,
	Beijing, China}
\email{yc-zhu@mail.tsinghua.edu.cn}

\date{}

\begin{document}
\maketitle

\section{\bf Introduction  }\label{section1}

The coinvariants and conformal blocks for affine Lie algebras was introduced in \cite{TUY}, where their geometric properties, including connections and factorization, were studied. The notions was subsequently extended to the more general frameworks of vertex operator algebras in \cite{FBZ} and chiral algebras in \cite{BD}. For a family of pointed curves with a module over a vertex operator algebra assigned to each marked point, one can construct a sheaf of coinvariants over the base, equipped, under suitable assumptions, with a projectively flat connection arising from the Virasoro field over the smooth locus. Finite-dimensionality and factorization results for coinvariants over the projective line were established in \cite{NT}, while finiteness results for conformal blocks over compact Riemann surfaces were proved in \cite{AN} under the $C_2$ condition. More recently, the geometric properties of sheaves of coinvariants over moduli spaces of stable pointed curves, including their connections and vector bundle structures, were investigated in \cite{DGT21}, while important factorization results were established in \cite{DGT24}. These results are obtained under various the representation-theoretic assumptions on the vertex operator algebras and their modules, with rationality and the $C_2$-cofiniteness condition playing important roles. We refer to these works and the references therein for further details. 

The notion of quasi-lisse vertex operator algebras, introduced in \cite{AK}, provides a broader framework than the $C_2$-cofiniteness condition and includes many interesting examples beyond the rational and $C_2$-cofinite setting. Motivated by this development, in the present work we study coinvariants of vertex operator algebras over the Weierstrass family of cubic curves, including their degenerations to nodal and cuspidal cubic curves. 

The Weierstrass family is 
\[
{\mathcal C}:\quad Y^2Z = 4X^3 - g_2 X Z^2 - g_3 Z^3,
\]
parametrized by $(g_2,g_3)\in\C^2$ and equipped with the zero section
$s\colon (g_2,g_3)\mapsto (0,1,0)$.  Under the identification
\[
g_2 = \frac{4\pi^4}{3}E_4, \qquad g_3 = \frac{8\pi^6}{27}E_6,
\]
the base is ${\rm Spec} \C[E_4,E_6]$, where $E_4,E_6$ are the Eisenstein
series and $\C[E_4,E_6]$ is the full ring of modular forms for
${\rm SL}_2(\Z)$.  The discriminant $\Delta = E_4^3 - E_6^2$ governs the
geometry of the fibers: when $\Delta \ne 0 $, the fiber is a smooth elliptic curve,  when $\Delta = 0$ and
$(E_4,E_6)\neq(0,0)$, the fiber is an irreducible nodal cubic curve, and the cuspidal cubic curve $y^2=4x^3$  when  $ (E_4 , E_6 )= (0 , 0 )$.

\medskip

Let $V$ be a vertex operator algebra with Virasoro element $\omega$.
Attaching $V$ at the zero section and taking coinvariants over the
Weierstrass family gives a module over $\C[E_4,E_6]$, which we call
the module of elliptic coinvariants:
\[
C_{\rm ell} (V) := \C[E_4,E_6]V \big/ \mathcal L(V)V,
\]
where $\mathcal L(V)$ is the $\C[E_4,E_6]$-span of the residue operators
$a_{(0)}$ and $\Res_z\bigl(Y(a,z)\wp_2(z,\tau)\bigr)$, $a\in V$, with
$\wp_2$ the Weierstrass $\wp$-function.  This is the space of
coinvariants of the Weierstrass family with $V$ assigned at $s$.  It carries two canonical operators: The weight operator assigns the total weight, combining the modular form weight and the conformal weight,
 and an operator
$D$ whose restriction to $\C[E_4,E_6]$ is the Ramanujan--Serre
derivative
\[
D = -\tfrac13 E_6\,\partial_{E_4} - \tfrac14 E_4^2\,\partial_{E_6}.
\]
The two satisfy $[W,D]=2D$, so that $C_{\rm ell}(V)$ is a 
$\C[E_4,E_6]\la D,W \ra$-module.  
We compute the action of $\mathcal L(V)$ through the Laurent expansions of the
$\wp_k$ and the quasi-periodicity of $\wp_1$, rather than through
algebro-geometric descent.

\medskip

Assume that $V$ is quasi-lisse, which implies that the subspace
\[
\operatorname{Span}( a_{(0)}b,\ a_{(-2)}b\mid a,b\in V ) 
\]
has finite codimension in $V$, a condition weaker than $C_2$-cofiniteness. Under this assumption, $C_{\rm ell}(V)$ is a finitely generated $\mathbb C[E_4,E_6]$-module and hence defines a coherent sheaf $C_{\rm ell}(V)^{\sim}$ over $\mathbb C^2$. The operators $D$ and $W$ induce a flat connection on the restriction of this sheaf to the complement of the discriminant locus $\Delta = 0 $.

\medskip

The torsion of $C_{\rm ell}(V)$ as an $\C[E_4,E_6]$-module ,  $C_{\rm ell}(V)_{\rm tor}$,
is also a $\C[E_4,E_6]\la D,W \ra$-module ( Lemma \ref{lemma4.1}). When $V$ is quasi-lisse, we show that 
the annihilator ideal of the torsion submodule $C_{\rm ell}(V)_{\mathrm{tor}}$
contains a power of  $\Delta \in \C[E_4,E_6]$
(Corollary \ref{cor4.3}).  Equivalently, the torsion is supported
 on the discriminant locus, i.e.\ on the union of the nodal and
cuspidal fibers, and away from that locus $C_{\rm ell}(V)$ is torsion-free.  When $V$ is rational and satisfies the $C_2$-condition,  
$C_{\rm ell}(V)_{\mathrm{tor}}$ is finite dimensional and $C_{\rm ell}(V)_{\mathrm{tor}}^{\sim}$ is supported at $( 0 , 0 )$ (Theorem \ref{theorem5.3}).

\

It is easy to  identify the fibers of $C_{\rm ell}(V)$ at the singular members of the
family.  They are related the invariants the associative algebra $A(V)$ and the Poisson algebra $R_V$ associated to $V$ introduced in \cite{Zhu}.
 The
fiber over a nodal cubic curve is naturally isomorphic to the cocenter of the
associative algebra $A(V)$, and that the fiber over the cuspidal cubic
is naturally isomorphic to the cocenter of the Poisson algebra $R(V)$
(Proposition \ref{prop5.2}).  These isomorphisms connect the geometry of
the Weierstrass family to the two classical invariants $A(V)$ and $R_V$ of $V$.

\

When $V$ is quasi-lisse,  the restriction $C_{\rm ell} ( V )^{\sim } $ to $\Delta \ne 0 $ is a holomorphic vector bundle with a flat connection, we get a
 monodromy representation of $ \pi_1 ( \Delta \ne 0 ) $ which is the braid group $B_3$, we prove that this representation factors through ${\rm SL}_2 ( \Z)$ (Theorem \ref{theorem4.4}).
On other hand side, it is known that $ Tr_M o( a) q^{ L_0 -\frac c {24} } $ converges on $0 < |q | <1 $ \cite{AK}. To study the modular invariance of these trace functions, we introduce the $ \C [ E_4 , E_6] \la D , W \ra$-module $\mathcal H$  formed by   
holomorphic functions on the upper half plane $\mathbb H$ (Section 5.3).  We show there is a natural action of ${\rm SL}_2 ( \Z )$ on 
 $\mathcal H$ as $ \C [ E_4 , E_6] \la D , W \ra$-module automorphisms (Proposition \ref{prop5.5}). 
We interpret the trace map as a homomorphism from
$C_{\rm ell}(V)$ to $\mathcal H$ of $ \C [ E_4 , E_6] \la D , W \ra$-modules  (Proposition \ref{prop5.4}). By the results in \cite{Zhu}, for rational $V$ that satisfies the $C_2$-condition, ${\rm Hom} ( C_{\rm ell}(V) ,  \mathcal H     )$
  has a basis labelled by the irreducible $V$-modules, the trace functions for each irreducible $V$-module gives a base element.  We show that 
  the monodromy representation of ${\rm SL}_2 ( \Z ) $ is the dual to the monodromy representation of ${\rm SL}_2 ( \Z ) $ (Theorem \ref{theorem5.6}).
  
 \
 
 This paper is organized as follows. Section 2 recalls the necessary background on modular forms, elliptic functions, and the Weierstrass family of cubic curves. In particular, we introduce the weight derivation $W$ and  the Ramanujan-Serre derivation $D$ on the space of modular forms. 
  We extend these derivations to the elliptic function algebra $\mathcal E$, which serves as the coordinate ring of the Weierstrass family with the zero section removed. In Section 3, we define the module of elliptic coinvariants $C_{\rm ell} (V)$
  for a vertex operator algebra $V$ and establish its basic properties, including the actions of 
$W$ and $D$.
  Section 4 is devoted to the study of $C_{\rm ell} (V)$
 under the quasi-lisse assumption: we prove that $C_{\rm ell} (V)$
  is a finitely generated $\C[ E_4 , E_6 ]$-module, describe its torsion submodule and its support on the discriminant locus, and show that the monodromy representation of the flat connection on the smooth locus factors through ${\rm SL}_2 ( \Z )$.
  In Section 5, we introduce $ \C [ E_4 , E_6] \la D , W \ra$-module $\mathcal H$ 
 of holomorphic functions on the upper half plane, and relate the trace functions to homomorphisms of 
 $ \C [ E_4 , E_6] \la D , W \ra$-modules. 
 
\

{\bf Acknowledgements}. This work is supported by the National Key R\&D Program of
China (Grant No. 2025YFA1017400).

\

\section{\bf Modular Forms, Elliptic Functions and Weirstrass Family of Cubic Curves }\label{section2}

 In this section, we recall some formulas concerning modular forms and elliptic functions, and then we introduce the Weierstrass family of cubic curves.
 
\subsection{Modular Forms of ${\rm SL}_2 ( \Z )$ and Ramanujan-Serre Derivative}
 
 For $k \geq 2 $, the Eisenstein series $G_{2k} ( \tau ) $ is 
 \begin{equation}
  G_{2k} ( \tau ) = \sum_{  ( m , n ) \in {\Z}^2 , ( m , n ) \ne ( 0 , 0 )}  \frac 1 { ( m \tau + n )^{2k} }
\end{equation}
  $G_{2k}$ is a modular form of weight $2k$ for ${\rm SL}_2 ( \Z )$.
 We set
 \begin{equation}
  E_{2k}(\tau ) = \frac 1 { 2\zeta (2k)} G_{2k} (\tau ) = 1 +\frac { ( 2 \pi i)^{2k} } { ( 2k-1) ! \zeta (2 k ) }
   \sum_{ n=1}^\infty   \sigma_{ 2k-1} ( n ) q^n.
 \end{equation}
 The  ring of modular forms for ${\rm SL}_2 ({\Bbb Z})$
 is the polynomial ring $ \C [ E_4 (\tau ) , E_6 (\tau )]$.
 
We need two fundamental derivations of the ring $ \C [ E_4 (\tau ) , E_6 (\tau )]$. The first derivation $W$ gives the weight of the 
 modular forms, so $ W E_4 = 4 E_4 ,  W E_6 = 6 E_6 $, it can written as
 \begin{equation}
      W = 4 E_4 \partial_{E_4} + 6 E_6 \partial_{E_6 }  . 
 \end{equation} 
 
 The second derivation we use is the Ramanujan-Serre derivative $D$. To give an expression of $D$ in terms of $\frac d { d \tau }$, 
 we need the quasi-modular forms $G_2 (\tau ) $ or $E_2 (\tau )$. They are defined as
 \[  G_2 (\tau ) = \frac { \pi^2 } 3 - 8 \pi^2  \sum_{  n =1}^\infty \sigma (n ) q^{n } \]
 \[ E_2 (\tau ) = \frac { 3 } { \pi^2 } G_2 ( \tau ) = 1 - 24  \sum_{  n =1}^\infty \sigma (n ) q^{n } .\]
 $G_2(\tau ) $ satisfies
 \begin{equation}\label{G2}
  G_2 ( \tau +1 ) = G_2 (\tau ) , \; \; \;    \tau^{-2}  G_2 ( -\frac 1 \tau ) = G_2 (\tau ) - \frac { 2 \pi i } { \tau } .
  \end{equation}
  Equivalently,
   \begin{equation}
 E_2 ( \tau +1 ) =  E_2 (\tau ) , \; \; \; \; \;  E_2 (-\frac 1 \tau )=\tau^2  E_2 (\tau ) - \frac { 6 \tau i}  {  \pi }
\end{equation}
  (\ref{G2}) implies that 
  \begin{equation} \label{2.6} 
    G_2 \left( \frac {a \tau + b } { c\tau + d }  \right) = ( c\tau + d )^2 G_2 (\tau ) - 2 \pi i c ( c\tau + d ) 
  \end{equation}
  for $ g \in SL_2 ( {\Bbb Z} )$.
The Ramunujan-Serre derivation $D$ of the graded ring $ \C [ E_4 (\tau ) , E_6 (\tau )]$ is defined as follows,
	if $ f (\tau )$ is a modular form of weight $ k $ for $SL_2 ( {\Bbb Z})$, then 
\begin{equation}
D f  =  \frac 1 { 2\pi i } f' ( \tau )- \frac k { 12} E_2  f .
 \end{equation} 
  is a modular form of weight $k+2 $. 
It is determined by the relations
\begin{equation}
 D ( f g ) = D( f ) g + f D g ,  \qquad D E_4 = - \frac 1 3 E_6 ,  \qquad  D E_6 = - \frac 1 4 E_4^2 .
 \end{equation}
So it can be written as
\[    D =  - \frac 1 3 E_6 \partial_{E_4} -   \frac 1 4 E_4^2   \partial_{E_6}. \]       
We have the relation 
\begin{equation}\label{1GD}
  [ W , D ]  = 2 D .
\end{equation}

\

\subsection{Elliptic Functions}

  Recall the 
 Weierstrass $\wp$-functions:  for fixed $\tau \in {\Bbb H}$, let $L_\tau $ denote the lattice
 $ L_\tau = {\Bbb Z} \tau + {\Bbb Z} $, $L_\tau' = L_\tau -\{0\} $,
 \begin{eqnarray*}
 && \wp_1 ( z , \tau ) = \frac 1 z + \sum_{ v \in L_\tau'} \left( \frac 1 { z - v } + \frac 1 v + \frac z { v^2 } \right) \\
 && \wp_2 ( z , \tau ) = \frac 1 {z^2} + \sum_{ v \in L_\tau'} \left( \frac 1 { ( z - v )^2 } - \frac 1  { v^2 }\right)  \\
 && \wp_{k+1} (z , \tau ) = - \frac 1 k \frac {d} {dz} \wp_k ( z , \tau ) .
 \end{eqnarray*}
  For $k \geq 2 $, $\wp_k (z , \tau )$ are periodic functions of $z$ with periods lattice $ L_\tau $.

 For fixes $\tau $, the functions $1$,  $\wp_k (z , \tau )$ ($k=2 , 3 , \dots $) form a basis of the meromorphic functions on the elliptic curve $ {\Bbb C} / L_\tau $
  with possible poles in $\{0\}$.

 The function  $\wp_1 ( z , \tau )$ has the properties:
 \begin{equation}\label{wp1}
  \wp_1 ( z +1  , \tau )= \wp_1 ( z , \tau ) + G_2 (\tau ) , \; \; \;
 \wp_1 ( z +\tau   , \tau )= \wp_1 ( z , \tau ) + G_2 (\tau )\tau - 2\pi i
   \end{equation}
 We have the Laurent expansion of $\wp_k$ for $k\geq 1 $ at $z=0$:
 \begin{equation}
  \wp_k ( z  , \tau )= \frac 1 { z^k} + (-1)^k \sum_{n=1}^\infty \binom{2n+1}{k-1} G_{2n+2} (\tau ) z^{2 n+2 - k }
  \end{equation}

 We introduce the {\it ellitic function algebra} 
\begin{equation}
\mathcal E=\C [ E_4 , E_6 ]\cdot 1\oplus\bigoplus_{n\geq0}\C [ E_4 , E_6 ]\cdot\partial_z^n\wp(z,\tau).
\end{equation}
The algebra $\mathcal E$ is the coordinate ring of the Weierstrass family of cubic curves with the infinity removed, see Section 2.3.  

\

We want to extend the derivations $ W $ and $D$ to $ \mathcal E$.	
	The ring 
$ {\mathcal E} $ has a grading   given by
\begin{equation}
{\rm deg}  \, E_4 = 4 , \; \; {\rm deg} \, E_6 = 6, \; \; {\rm deg} \partial_z^n \wp_2 ( z , \tau ) = 2 + n  
\end{equation} 
 So we have
\[  {\mathcal E} = \oplus_{ n = 0 }^\infty  {\mathcal E}_n .     \]
It is easy to see that 
$ f ( z , \tau ) \in  {\mathcal E}$ is in ${\mathcal E}_n$ if and only if 
\[     f \left( \frac z{ c \tau + d } , \frac { a \tau + b} { c \tau + d }\right) = ( c\tau + d)^n f ( z , \tau ) \qquad {\rm for} \; \;    \left( \begin{matrix} a & b \\ c & d \end{matrix} \right) \in {\rm SL}_2 ( \Z) .\]
This grading gives the derivation $W$: 
\[   W f = n f \qquad {\rm for}  \; f \in {\mathcal E}_n , \]
which extends the derivation on  $ \C [ E_4 , E_6]$ denoted by the same symbol. 

\

To extend $D$ to $\mathcal E$, we first introduce a differential operator $X$ on the meromorphic functions on $ \C \times {\mathbb H} $. 
\[  X = ( 2\pi i )^{-1}\partial_\tau  + ( 4 \pi^2)^{-1} \wp_1 ( z , \tau ) \partial_z - \frac 1 { 12 } E_2 (\tau ) z \partial_z .\]

\

\begin{lemma}\label{lemma2.1}
	Let $f ( z , \tau ) $ be a meromorphic functions on ${\Bbb C} \times {\Bbb H}$ satisfying the elliptic property
	\begin{equation}
	f ( z + 1, \tau ) = f ( z , \tau ) , \; \; f ( z + \tau , \tau )= f ( z , \tau ).
	\end{equation}
	Then $Xf$ satisfies  the similar properties:
	\[  (X f )( z + 1, \tau ) = ( Xf ) ( z , \tau ) , \; \; (X f ) ( z + \tau , \tau ) = (X f ) ( z, \tau ), \]
\end{lemma}

\begin{proof} 
	Let 
\begin{equation}\label{2g}
 \tilde X = (2\pi i )^2 X = 
 2\pi i \partial_\tau  - \bigl( \wp_1 ( z , \tau ) - G_2 (\tau ) z \bigr)  \partial_z =   2\pi i \partial_\tau  - A ( z , \tau )  \partial_z  ,\end{equation}
where
	\[
	A(z,\tau)=\wp_1(z,\tau)-G_2(\tau)z.
	\]
	By the quasi-periodicity formulas for $\wp_1$, we have  
	\begin{equation}\label{eq:Aperiods}
	A(z+1,\tau)=A(z,\tau),\qquad
	A(z+\tau,\tau)=A(z,\tau)-2\pi i.
	\end{equation}
	Differentiating the periodicity identities for $f$ with respect to $z$ gives
	\begin{equation}\label{eq:zder}
	\partial_z f(z+1,\tau)=\partial_z f(z,\tau),\qquad
	\partial_z f(z+\tau,\tau)=\partial_z f(z,\tau).
	\end{equation}
	Differentiating $f(z+1,\tau)=f(z,\tau)$ with respect to $\tau$ gives
	\begin{equation}\label{eq:tau1}
	\partial_\tau f(z+1,\tau)=\partial_\tau f(z,\tau).
	\end{equation}
	Equations \eqref{eq:Aperiods}--\eqref{eq:tau1} imply
	$(\tilde Xf)(z+1,\tau)=(\tilde Xf)(z,\tau)$.
	
	For the second period, differentiate $f(z+\tau,\tau)=f(z,\tau)$
	with respect to $\tau$, taking account of the dependence of $z+\tau$
	on $\tau$. The chain rule yields
	\[
	\partial_z f(z+\tau,\tau)+\partial_\tau f(z+\tau,\tau)
	=\partial_\tau f(z,\tau),
	\]
	so that
	\begin{equation}\label{eq:tau2}
	\partial_\tau f(z+\tau,\tau)
	=\partial_\tau f(z,\tau)-\partial_z f(z,\tau).
	\end{equation}
	Using \eqref{eq:Aperiods}, \eqref{eq:zder}, and \eqref{eq:tau2}, we obtain
	\begin{align*}
	(\tilde Xf)(z+\tau,\tau)
	&=2\pi i\,\partial_\tau f(z+\tau,\tau)
	-A(z+\tau,\tau)\partial_z f(z+\tau,\tau)\\
	&=2\pi i\bigl(\partial_\tau f(z,\tau)-\partial_z f(z,\tau)\bigr)
	-\bigl(A(z,\tau)-2\pi i\bigr)\partial_z f(z,\tau)\\
	&=2\pi i\,\partial_\tau f(z,\tau)-A(z,\tau)\partial_z f(z,\tau)\\
	&=(\tilde Xf)(z,\tau).
	\end{align*}
	This proves both periodicity identities.
\end{proof}

\

\begin{lemma}
	 We have
\begin{equation}\label{2p2}
 X \wp_2 ( z , \tau ) =  - ( 2 \pi^2)^{-1} \wp_4 ( z , \tau )  + (2 \pi^2)^{-1} G_2 ( \tau ) \wp_2 ( z , \tau ) + 5 ( 2 \pi^2)^{-1} G_4 (\tau )
\end{equation}
and for every $ n \geq 0 $,
 \begin{equation}\label{2.19}
 \left( X  - \frac 1 {12} ( n +2 ) E_2 (\tau ) \right) \partial_z^n\wp_2 ( z , \tau ) \in \mathcal E
\end{equation}
\end{lemma}

\begin{proof}  Recall operator $\tilde X = - 4 \pi^2 X$, see (\ref{2g})
	\begin{eqnarray*}
		(\tilde X \wp_2 ) ( z , \tau ) &=&  \left( 2 \pi i \partial_\tau - ( \wp_1 (  z, \tau ) - G_2 (\tau ) z )\partial_z \right) \wp_2 ( z , \tau )\\
		&=&   \left( 2 \pi i \partial_\tau - (z^{-1} - \sum_{ n=0}^\infty  G_{2n+2} (\tau ) z^{2n+1} )\partial_z \right) \wp_2 ( z , \tau )\\
		&=&    2 \pi i \partial_\tau \wp_2 ( z , \tau )  + 2 (z^{-1} - \sum_{ n=0}^\infty  G_{2n+2} (\tau ) z^{2n+1} ) \wp_3 ( z , \tau )\\
		&=&  2 \pi i \partial_\tau \wp_2 ( z , \tau )  + 2 (z^{-1} - \sum_{ n=0}^\infty  G_{2n+2} (\tau ) z^{2n+1} ) ( z^{-3} -\sum_{ n=1}^\infty {{2n+1} \choose 2 } G_{2n+2}z^{2n-1} )\\
		&=&     2  z^{-4} -2 G_2 (\tau )z^{-2} - 8 G_4 ( \tau ) + {\rm sums \; of \; positive \; powers \; of \; }  z \\
		&=&  2 \wp_4 ( z , \tau ) -2 G_4 ( \tau ) - 2 G_2 ( \tau ) \wp_2 ( z , \tau )  - 8 G_4 (\tau ) +  {\rm sums \; of \; positive \; powers \; of \; }  z 
	\end{eqnarray*}
By Lemma \ref{lemma2.1},  as a function of $z$, $  (\tilde X \wp_2 ) ( z , \tau ) $ is doubly periodic with periods $1 $ and $\tau $, so
	\[   (\tilde X \wp_2 ) ( z , \tau ) - \left( 2 \wp_4 ( z , \tau )  -2  G_2 ( \tau ) \wp_2 ( z , \tau ) - 10 G_4 (\tau )       \right) \]
it is doubly periodic with periods $1$ and $\tau $, it also an entire function of $z$, and its value at $0$ is 0, so the function is $0$. This proves
	\begin{equation}
\tilde	X \wp_2 ( z , \tau ) =  2 \wp_4 ( z , \tau )  - 2 G_2 ( \tau ) \wp_2 ( z , \tau ) - 10 G_4 (\tau ). 
	\end{equation}
	which is equivalent to (\ref{2p2}). To prove (\ref{2.19}), we use induction on $n$. The case $n=0 $ follows from (\ref{2p2}).  
 The proof of the induction step is proved	using the identity 
	\[        [ X , \partial_z ]= - ( 4 \pi^2)^{-1} \wp_2 ( z , \tau ) \partial_z + \frac 1 { 12 } E_2 ( \tau )\partial_z .\]
 \end{proof}

The lemma shows that $ X$ does not preserve $\mathcal{E}$, we modify it as 
 \begin{equation}
 D = X - \frac 1 {12} E_2 (\tau) W  = ( 2\pi i )^{-1}\partial_\tau  + ( 4 \pi^2)^{-1} \wp_1 ( z , \tau ) \partial_z - \frac 1 { 12 } E_2 (\tau ) z \partial_z - \frac 1 { 12 } E_2 (\tau ) W.
 \end{equation}

\

\begin{lemma}\label{lemma2.3}
 The operator $D$  a derivation of $ \mathcal E$. Its restriction on $\C [ E_4 , E_6 ]$ is the Ramnujan-Serre derivative. 
 And we have
 \[  [ W , D ] =2 D , \; \;  [ W , \partial_z ] = ( 4 \pi^2)^{-1}\wp_{2 } ( z , \tau ) \partial_z   \]
\end{lemma}

\begin{proof} for modular form $ f (\tau )$ of weight $k$, since $\partial_z f = 0 $,
	 $ Df = ( 2 \pi i)^{-1} f' - \frac k {12} E_2 f $ is the Ramanujan-Serre derivative of $f$.
	By (\ref{2.19}),  $ D \partial_z^n \wp_2 ( z , \tau ) \in \mathcal E$.  Since  $ \partial_z^n \wp_2 ( z , \tau ) $ is a basis 
	 of $ \C [ E_4 , E_6 ]$, $D$ preserves $\mathcal E$.  The other claims in the lemma are trivial.
	 \end{proof}

\

\subsection{Weierstrass family of Cubic Curves}

We consider the family of cubic curves parametrized by  $ (g_2 , g_3 ) \in \C^2 $:
\[ {\mathcal C}: \; \;     Y^2 Z = 4 X^3 - g_2 X Z^2 - g_3  Z^3 ,\]
We have the section 
\[ s:  \C^2 \to \mathcal C , \qquad ( g_2 , g_3) \mapsto  O  = (0 , 1 ,  0 ) \]
 We call the family  $\mathcal C $ with the section $s$ 
 the Weierstrass family of cubic curves. See \cite{Sil} for further details. 

 To connect it with modular forms,
we set 
\begin{equation}
g_2=60G_4=\frac{4\pi^4}{3}E_4,\qquad
g_3=140G_6=\frac{8\pi^6}{27}E_6.
\end{equation}
The $\C [ E_4 , E_6 ]$-algebra ${\mathcal E}$ is 
isomorphic to the quotient algebra 
\[	 	\C [ E_4 , E_6 ] [ X , Y ]/(Y^2-4X^3+g_2X+g_3) \]
under the map $ X \mapsto \wp_2 (z , \tau ), Y \mapsto \partial_z \wp_2 (z , \tau ) $. The later is 
the coordinate ring of $ \mathcal C - s $. 
The operators $W$ and $D$ gives the vector fields on  $ \mathcal C - s $

We write $\Delta (\tau ) =E_4 (\tau )^3-E_6(\tau )^2$, which spans the cusp form for $SL_2 (\Z)$ with weight $12$. The properties of the fibers
\[  \pi :   {\mathcal C} \to \C^2 = {\rm Spec} \, \C [ E_4 , E_6 ] \]
 are as follows:
\newline (1) If $\Delta (\tau )\ne0$, the fiber $ \pi^{-1} ( E_4 ( \tau ) , E_6 ( \tau )) $  is an elliptic curve with distinguished point $O$.
\newline (2) If $\Delta (\tau)=0$ but $(E_4 (\tau ),E_6 (\tau ) )\ne(0,0)$, the fiber is an irreducible nodal cubic curve.
\newline (3) At $(E_4,E_6)=(0,0)$, the fiber is the cuspidal cubic curve $y^2=4x^3$.

In particular, $\pi$ is a smooth family over  
\[
S^{\mathrm{sm}}=\operatorname{Spec}\C [ E_4 , E_6 ][ \Delta^{-1}],
\]
but not over all of $\C^2$. The nodal fibers are stable as one-pointed genus-one curves; the cuspidal fiber is not stable.

The complex analytic fiber for the case $\Delta \ne 0 $ can described as follows. For $\tau$ in the upper half-plane $\mathbb H$, let
\[
\Lambda_\tau=\mathbb Z+\tau\mathbb Z,
\qquad E_\tau=\mathbb C/\Lambda_\tau.
\]
The Weierstrass elliptic function satisfies
\[
\wp'(z,  \tau)^2
=4\wp(z , \tau)^3
-g_2(\tau)\wp(z,  \tau)-g_3(\tau).
\]
Consequently, the fiber over the point is isomorphic to $E_\tau $, the point $ O $ is $ 0 \in E_\tau $.

\

\section{\bf  Module of Elliptic Coinvariant for Vertex Operator Algebras }\label{section3}

\subsection{Definition of Elliptic Coinvariant Module}\label{section3.1}

Let $V $ be a vertex operator algebra with Virasoro element $\omega $, so
\[ Y ( \omega , z ) = \sum_{ n \in {\Bbb Z} } \omega_{(n)} z^{-n-1} = \sum_{ n \in {\Bbb Z} } L_n z^{-n-2}. \]
Suppose $ V = \oplus_{n=0}^\infty V_n $ with $ V_n = \{ v \, | \, L_0 v = n v \}$.
We consider the space 
$   \C [ E_4 , E_6 ] \otimes V =   \C [ E_4 , E_6 ]  V  $.
We define a $\C [ E_4 , E_6 ]$-linear map
\begin{equation}\label{3.0}
   V \otimes \mathcal E \to {\rm End}_{ \C [ E_4 , E_6 ] }  (   \C [ E_4 , E_6 ] V ) , \qquad  a \otimes f ( z , \tau ) \mapsto
{\rm Res}_z \left( Y ( a, z ) f (z , \tau ) \right) .
\end{equation}
By the formula
\[     [   {\rm Res} \left( Y (a , z) f( z ) \right) , \;     {\rm Res} \left( Y (b , z) g ( z ) \right)       ]
= \sum_{ i = 0}^\infty \frac 1 { i !} {\rm Res} \left(  Y ( a_{(i)} b , z )  f^{(i)} (z) g ( z )  \right) \]
and the fact that $ \mathcal E$ is an algebra closed under $\partial_z $,  the span of the image of (\ref{3.0}) is a Lie algebra, we denote it by
 $ {\mathcal L} ( V) $.  Since $ \mathcal E$ is generated by $1$, $\wp_k ( z, \tau )$ ($k \geq 2 $) as a $ \C [ E_4 , E_6 ]$-module,
  $ {\mathcal L} ( V) $ is 
 the $ \C [ E_4 , E_6 ]$-span of 
\[   {\rm Res} \left( Y ( a , z ) 1 \right)  = a_{(0)}  \]
and 
\begin{equation}\label{3.1}
    {\rm Res} \left( Y ( a , z ) \wp_k (z , \tau )  \right)
     = a_{(-k)} + (-1)^k \sum_{n=1}^\infty \binom{2n+1}{k-1} G_{2n+2} (\tau ) a_{ (2 n+2 - k ) } , \qquad k \geq 2 
\end{equation}     
Since 
\[ \left( \frac d { d z } \right)^k \wp_2 ( z , \tau  ) =  (-1)^k ( k +1 ) !   \wp_{2 + k } ( z , \tau ) \]
we have
\begin{align*}
 &  (-1)^k ( k +1 ) ! {\rm Res} \left( Y (a , z)   \wp_{2 + k } ( z , \tau )  \right) \\
  &    =
    {\rm Res}\left( Y (a , z) \left(  \frac d { d z } \right)^k \wp_2 ( z , \tau  )  \right) \\
  &  =(-1)^k    {\rm Res}\left( Y ( L_{-1}^k a , z) \wp_2 ( z , \tau  )  \right) 
\end{align*}
From this formula, we see that 
${\mathcal L} ( V)$ is the $ \C [ E_4 , E_6 ]$-span of $a_{(0)} $ and $  {\rm Res} \left( Y ( a , z ) \wp_2 (z , \tau )  \right)$, $  a\in V $.

\

\begin{defn}
	The module of elliptic coinvariant for $V$  is the quotient module  
	\[  C_{\rm ell} ( V) \deff \C [ E_4, E_6 ]  V /      {\mathcal L} ( V )   V  . \]
\end{defn}

 The space $C_{\rm ell} ( V)$ is the space of coinvariants associated the Weierstrass cubic family with $V$ assigned at the section $s$.

\

We introduce a grading on $  \C [ E_4 , E_6 ]  V  $ by
\[  {\rm deg}  (  f ( \tau ) \otimes  a ) = k + n , \]
where $ L_0 a = n a $ and $ k $ is the weight of the modular form $ f $. For example, ${\rm deg} ( E_6 \otimes \omega ) =8 $.
We define the corresponding degree operator $ W $ by 
\[   W  \,  v =  {\rm deg} ( v)  v   \] 
for homogeneous element $v$ and extend $W$ to whole the space by linearity. 
For $a \in V_m ,  b \in V_n $,
$ a_{(0)} b \in V_{ m + n-1 }$.  Every term in summation 
\[   {\rm Res} \left( Y ( a , z ) \wp_2 (z , \tau )  \right) b  = a_{(-2)} b +  \sum_{n=1}^\infty (2n+1 ) G_{2n+2} (\tau ) a_{ (2 n ) } b \]
has degree $ m + n +1 $, as $ G_{2n+2}$ has degree $2n +2 $.  The subspace   ${\mathcal L} ( V )   V$ is a graded subspace, therefore the grading
 on  $\C [ E_4, E_6 ]  V$ descends to $  C_{\rm ell} ( V) $. We have
 \[      C_{\rm ell} ( V) = \oplus_{ n = 0}^\infty    C_{\rm ell} ( V)_n                   \]
 where  $C_{\rm ell} ( V)_n $ is the subspace space formed by the elements of degree $n$.  
    We denote the corresponding degree operator again by  $W$, so 
 $ W |_{ C_{\rm ell} ( V)_n } =  n \, {\rm Id } $. 
 
 \begin{lemma}\label{lemma3w}
 	For $ a\in V_m $, that is, $L_0a = m a $, $ f ( z , \tau ) \in \mathcal E_n $, we have
 \begin{equation}\label{3w}
   [ W ,  {\rm Res}_{z=0} \left(  Y (a , z ) f ( z , \tau ) \right) ] = 
  ( m + n -1 )  {\rm Res}_{z=0} \left(  Y (a , z ) f ( z , \tau ) \right)   \end{equation}
 \end{lemma}
 
 \begin{proof}
 	It is enough to prove (\ref{3w}) for the case $ f =1 $ and $ f =\wp_k $, the result follows from the formula (3.2).
 \end{proof}

\

\subsection{The Operator $D$}\label{section3.2}

We define an operator $ D $ on $ \C [ E_4 , E_6 ]  V$ by
\begin{align}
  D  & = ( 2 \pi i )^{-1} \partial_\tau + \frac 1 { 12 } E_2 (\tau ) L_0 - \frac 1 {12}  E_2 (\tau ) W - ( 4 \pi^2 )^{-1}
    {\rm Res}_z \left(     Y (\omega , z )\wp_1 ( z , \tau ) \right)
   \nonumber \\
     & =    (2 \pi i)^{-1} \partial_\tau  +  \frac 1 { 12} E_2 (\tau )L_0 - \frac 1 {12} E_2(\tau ) W  - ( 4\pi^2)^{-1} L_{-2} +  
     ( 4\pi^2)^{-1} \sum_{ n = 1 }^\infty G_{2n+2} L_{2n}  
\end{align}
  The operator $D$ can also be written as
  \begin{equation}
  D =  ( 2 \pi i )^{-1} \partial_\tau  - \frac 1 {12}  E_2 (\tau ) W - ( 4 \pi^2 )^{-1}
  {\rm Res}_z \left(     Y (\omega , z )( \wp_1 ( z , \tau ) - G_2 (\tau ) z )\right)
  \end{equation}

Although $G_2(\tau ) $ appears in the expression of $D$, we show that $ D$ acts on the space $ \C [ E_4 , E_6 ]  V$. For $ f (\tau ) v \in \C [ E_4 , E_6 ]  V$,
where $ f $ is a modular form of weight $k$,  $v \in V_m $, we have $ {W} ( f v ) = ( k + m ) f v $ by the definition of $W$, 
\[  D ( f (\tau ) v  ) = ( (2 \pi i )^{-1} f' ( \tau ) - \frac k {12} E_2 (\tau ) f ) v  - ( 4\pi^2)^{-1} f (\tau ) L_{-2} v
 + (4\pi^2)^{-1} f(\tau )\sum_{ n = 1 }^\infty G_{2n+2} L_{2n}\]
since $  (2 \pi i )^{-1} f' ( \tau ) - \frac k {12} E_2 (\tau ) f  $ is a modular form of weight $k+2 $,  $ D  ( f (\tau ) v  )\in  \C [ E_4 , E_6 ]  V$.
This proves that $D$ is a $\C$-linear operator on $ \C [ E_4 , E_6 ]  V$.

\

\begin{prop}\label{prop3.3} 
	The operator $ D $ satsifies the following properties	
\newline (1) For a modular form $f ( \tau ) $, we denote $f$ the multiplication operator by $f$, then we have 
\[  [ D , f( \tau ) ] =  (D f ) (\tau ), \]
 where $ Df$ is the Ramanujan-Serre derivative of $f$.
\newline (2)  $D$ normalizes $ {\mathcal L} ( V )$, that is,  
	 \[  [ D , {\mathcal L} ( V )] \subset {\mathcal L} ( V ) \]
\newline (3) $D$ descends to a $\C$-linear operator   on  $ C_{\rm ell} ( V) $ and
 \[  D \, C_{\rm ell}( V)_n \subset C_{\rm ell}( V)_{n+2 }, \qquad   [ W , D] = 2 D .\]
\end{prop}

\

\begin{proof}
	To prove (1), for a modular form of weight $k$, we have 
	$ [ W, f (\tau )  ] = k $,
	\[  [ D , f(\tau ) ] = [ (2\pi i )^{-1} \partial_\tau - \frac 1 {12} E_2 (\tau) W  , f ( \tau )] = ( 2 \pi i )^{-1} f' (\tau ) -\frac k { 12 } E_2 (\tau ) f ( \tau ) = D f .\] 	
	To prove (2),  we need to prove, for an elliptic function $ f ( z , \tau )$ with ${\rm deg} \, f = k $, 
	\begin{equation}\label{3normal}
	 [ D ,  {\rm Res} \left( Y (a , z ) f (z  , \tau ) \right) ] \in {\mathcal L} ( V). 
	  	\end{equation}
	  	 We may assume  $ a \in V_m $ and ${\rm Deg} \, f = k $.  Let $ g ( z , \tau ) = \wp_1 ( z , \tau ) - G_2 (\tau ) z $,
	we have 
\begin{align*}
   & [ {\rm Res} \left( Y (\omega , z ) g ( z , \tau ) \right) ,     \; {\rm Res} \left( Y (a , z ) f (z  , \tau ) \right)  ]  \\
 &	=  \sum_{ j=0}^\infty \frac 1 { j!}  {\rm Res }  \left( Y ( g_{(j)} a , z ) g^{(j)} (z , \tau ) f ( z , \tau )   \right)  \\
 & = {\rm Res} \left( Y ( L_{-1} a , z ) \wp_1 (z , \tau ) f ( z , \tau ) \right) +  \sum_{ j=2}^\infty \frac 1 { j!}  {\rm Res }  \left( Y ( \omega_{(j)} a , z ) \wp_1^{(j)} (z , \tau ) f ( z , \tau )   \right) 
 \end{align*}
  We use  the formula
\begin{align*}	
& 	{\rm Res} \left( Y ( L_{-1} a  , z ) \wp_1 ( z , \tau   ) f( z , \tau ) \right) \\
&	= 	{\rm Res} \left( Y ' (  a  , z ) \wp_1 ( z , \tau   ) f( z , \tau ) \right)\\
&	= - {\rm Res} \left( Y  ( a  , w )  \wp_1 ( z , \tau   ) f'( z , \tau )\right) - {\rm Res} \left( Y  ( a  , w )  \wp_1' ( z , \tau   ) f(z , \tau ) \right)
	\end{align*}
Since $ \wp_1^{(j)} (z , \tau )$ is an elliptic function for $ j \geq 1 $, we have
\[     \sum_{ j=1}^\infty \frac 1 { j!}  {\rm Res }  \left( Y ( \omega_{(j)} a , z ) \wp_1^{(j)} (z , \tau ) f ( z , \tau )   \right) , \;
{\rm Res} \left( Y  ( a  , w )  \wp_1' ( z , \tau   ) f(z , \tau ) \right) \]
are in $ {\mathcal L} ( V) $, therefore we have
\begin{equation}\label{3mod}
[ {\rm Res} \left( Y (\omega , z ) \wp_1 ( z , \tau ) \right) ,     \;{\rm Res} \left( Y (a , z ) f (z  , \tau ) \right)  ] \equiv 
  - {\rm Res} \left( Y  ( a  , w )  \wp_1 ( z , \tau   ) f'( z , \tau )\right) \qquad   {\rm mod} \; {\mathcal L} ( V).
\end{equation}   

\begin{align*}
&   [ L_0 ,   {\rm Res} \left( Y (a , z ) f (z  , \tau ) \right)       ] \\
& =  {\rm Res} \left( Y (L_{-1}a , z )  z f (z  , \tau ) \right) + m  {\rm Res} \left( Y (a , z )   f (z  , \tau ) \right) \\
& =  - {\rm Res} \left(  Y (a , z )  z f' (z , \tau )  \right)   - {\rm Res} \left( Y (a , z )   f (z  , \tau ) \right) + m  {\rm Res} \left( Y (a , z )   f (z  , \tau ) \right) 
\end{align*}
So we have
\begin{equation}\label{3.2mod}
[ L_0 ,     \; {\rm Res} \left( Y (a , z ) f (z  , \tau ) \right)  ] \equiv  - {\rm Res} \left(  Y (a , z )  z f' (z , \tau )  \right) \qquad {\rm mod} \; \mathcal L ( V) 
\end{equation}
We also have
\[  [ \partial_\tau , {\rm Res } \left(Y (a , z ) f (z , \tau )  \right) ] = {\rm Res } \left(Y (a , z ) \partial_\tau f (z , \tau )  \right) \]
this together with Lemma \ref{lemma3w}, (\ref{3mod}) (\ref{3.2mod}) imply that
\begin{equation}
  [  D ,   {\rm Res} \left( Y (a , z ) f (z  , \tau ) \right) ] \equiv   {\rm Res} \left( Y (a , z ) D f (z  , \tau ) \right) \qquad {\rm mod} \; \mathcal L ( V) 
 \end{equation}
 by Lemma \ref{lemma2.3}, the right hand side is in $ \mathcal L ( V)$. This proves  (\ref{3normal}), therefore (2). 
 (3) follows from (2) directly. 
\end{proof}

\

\subsection{Algebra of Modular Differential Operators}\label{s3.3}

We introduce the algebra 
\begin{equation} 
 \C [ E_4 , E_6 ] \la D, W \ra = \oplus_{ m, n= 0}^\infty   \C [ E_4 , E_6 ] D^m W^n ,
\end{equation}
which has the algebra structure given by the commtation relation 
\[  D f - f D =  D(f) , \; \; \; W f - f W = W(f) ,  \; \; \;   [ W , D] = 2 D, \; \;     f \in \C [ E_4 , E_6 ] \]
We call this algebra the algebra of modular differential operators.

\

\begin{prop}
$ C_{\rm ell} ( V) $ is a module over  $\C [ E_4 , E_6 ] \la D , W \ra $, where $\C [ E_4 , E_6 ]$ acts as multiplication, $D$ and $W$ acts as the operators with the same symbols. 
If $V$ is quasi-lisse, then $ C_{\rm ell} ( V) $ is a finitely generated as a module over $\C [ E_4 , E_6 ]$. 
\end{prop}

\begin{proof}
	By Proposition \ref{prop3.3},  $ C_{\rm ell} ( V) $ is a module over  $\C [ E_4 , E_6 ] \la D , W \ra $.  If $V$ is quasi-lisse, by a result of 
	 Arakawa and Kawasetsu \cite{AK},  ${\rm Span} ( a_{(0)} b ,   a_{(-2)} b , a , b \in V) $ has finite codimension in $V$. Because the elements in 
	$ {\mathcal L} ( V )   V $ is spanned by $ a_{(0)} b$ and 
	\[  {\rm Res} \left(Y (a , z) \wp_2 ( z , \tau ) \right) = a_{(-2)} b + \sum_{ n =1}^\infty (2 n +1 ) G_{2n+2} (\tau ) a_{(2n)} b, \]
	their leading terms are $ a_{(0)} b $ and $a_{(-2)} b$,  $ C_{\rm ell} ( V) $ is a finitely generated as a module over $\C [ E_4 , E_6 ]$. 
\end{proof}

\

\section {Monodromy Representation of Coinvariants for Quasi-Lisse Vertex Operator Algebras }\label{section4}

In this section, we first study some general properties of modules over the algebra $\mathbb C[E_4,E_6]\langle D,W\rangle$. We show that the torsion submodule $M_{\rm tor}$ of such a module $M$ is stable under the actions of $D$ and $W$. If $M$ is finitely generated over $\mathbb C[E_4,E_6]$, we further prove that there exists an integer $m\geq 0$ such that $\Delta^m M_{\rm tor}=0$.

For a quasi-lisse vertex operator algebra $V$, the elliptic coinvariant module $C_{\mathrm{ell}}(V)$ defines a coherent sheaf $ {C_{\mathrm{ell}}(V)}^{\sim}$ on $\mathbb C^2$. We show that the operators $D$ and $W$ induce a flat connection on its restriction to the complement of the discriminant locus ${\Delta=0}$. Consequently, this restriction is a holomorphic vector bundle equipped with a flat connection. We then study the associated monodromy representation and prove that it factors through $\mathrm{SL}_2(\mathbb Z)$.

\subsection{\bf Structure of $\C [ E_4, E_6 ]$-finite  $\C [ E_4, E_6 ] \la  D ,  W \ra $-Modules }\label{section4.1}

We give a general discussion of the structure $\C [ E_4, E_6 ]$-finite  $\C [ E_4, E_6 ] \la  D ,  W \ra $-Modules.

 \

 \begin{lemma}\label{lemma4.1}
 	For a  $\C [ E_4, E_6 ] \la  D ,  W \ra $-Module $M$, let
 \[ M_{\rm tor} = \{   v \in M \, | \,  f v  = 0 \; \; {\rm for \; some \; } f \in \C [ E_4, E_6 ]  , f \ne 0 \} \]
 then $ M_{\rm tor} $ is a   $\C [ E_4, E_6 ] \la  D ,  W \ra $-submodule.
 \end{lemma}

\noindent {\it Proof. } If $ v \in M_{\rm tor } $, then $ f v = 0 $ for some $ f \ne 0 $, $ D f v = 0 $, $ D (f ) v + f D v = 0 $, multiply $f$ to this equation, we get 
$     f^2 D v  = 0 $. Similarly, $ f^2 W v = 0 $.
\hfill $\Box $

\

We identify $  \C [ E_4 , E_6 ] = \C [ x , y ] $ by $ E_4 \mapsto x $, $ E_6 \mapsto y $. Since 
 $W E_4 = 4 E_4 ,   W E_6= 6 E_6 $, $W$ is identified with 
\[ W  =  4 x \partial_x + 6 y \partial_y \]
The Ramanujan -Serre derivation  $D $ is identified with
\[ D = - \frac 1 3 y  \partial_x - \frac 1 2 x^2 \partial_y. \]
The transition matrix between $D , W $ and $\partial_x , \partial_y $ is
\[       \left( \begin{matrix}   - \frac 1 3 y & - \frac 12 x^2 \\  4 x & 6 y \end{matrix} \right) \]
its determinant is 
\[         2 ( x^3 - y^2 ) ,\]
which corresponds to $ 2 ( E_4^2 - E_6^2 ) = 2  \Delta $. Since
$ \partial_x $ and $\partial_y $ can be expressed as a $\C [ x , y ] [ ( x^3 - y^2 )^{-1}]$-linear combination of $W$ and $D$, therefore $M$  give a 
$D$-module $\widetilde M $ on $ \C^{2} - \{ y^2 - x^3 = 0 \} $.

\

\begin{prop}\label{prop4.2}
	Suppose $ M$ is a $\C [ E_4, E_6 ]$-finite  $\C [ E_4, E_6 ] \la  D ,  W \ra $-module. Then
\newline (1)	\[  M [\Delta^{-1}  ] = M \otimes_{\C [ E_4, E_6 ] }\C [ E_4, E_6 , \Delta^{-1}] \]
 is a projective module over $\C [ E_4, E_6 ] [ \Delta^{-1}]$ of finite rank. 
\newline (2)	There is a positive integer $m$ such that 
	\[  ( x^3 - y^2 )^m v = 0  \; \; {\rm for \; all \; }  v \in    M_{tor}. \] 
  \end{prop}

\begin{proof}  Let $ R =   \C [ x , y ] [ ( x^3 - y^2)^{-1}]$,  $ U = {\rm Spec} \, R =  \C^{2} - \{ y^2 - x^3 = 0 \} $.  
	$ {M [ ( x^3 -y^{2})^{-1}]}^{\sim}$ is  a coherent $D$-module over $U$. Since $U $ is smooth,  this sheaf is locally free.
	 Consequently, $M[\Delta^{-1}]$ is a finitely generated projective module over $R[\Delta^{-1}]$. 
In particular 
     $ M [\Delta^{-1}  ]_{\rm tor} = 0 $. 
 For (2), suppose $ v \in M_{tor}$,  then $ v \in M [\Delta^{-1}  ]_{tor} = 0 $, this implies that $ \Delta^n v = 0 $ for some $n\geq 1 $.
  Because $M_{tor}$ is finite over $\C [ E_4, E_6]$, we can find a uniform $m$ such that $( x^3 - y^2 )^m v = 0  $ for all $v\in M_{tor}$
\end{proof} 

\

It is easy to check that  
  \begin{equation}
     D ( E_4^3 - E_6^2 ) = 0 .
  \end{equation}
  
 \
 
 \begin{cor}\label{cor4.3}
   If $V$ is a quasi-lisse vertex algebra, then $  C_{\rm ell} (V) [ \Delta^{-1}] $ is a projective module over $\C [ E_4 , E_6 , \Delta^{-1}] $ of finite rank,
    $  C_{\rm ell} (V)_{tor}$ is annihilated by $\Delta^m $ for some $ m $.
 \end{cor}

  \
  
  \subsection{\bf Monodromy Representation of $C_{\rm ell} ( V) $ for Quasi-lisse VOA }\label{s4.2}
  
   We let 
   \[
   R =\C [ E_4 , E_6 ,  ( E_4^3-E_6^2 )^{-1} ] ,\qquad
   U={\rm Spec} R ,  \qquad    U ( \C )= \{ ( x , y ) \in \C^2 \, | \,   x^3 - y^2 \ne 0 \}.
   \]
 If $V$ is quasi-lisse, ${\mathcal C}_{\rm ell} ( V) $ a $ \C [ E_4 , E_6]$-finite $ \C [ E_4 , E_6 ] \la D, W \ra $-module. So we have 
 $ {\mathcal C}_{\rm ell} ( V) [ \Delta^{-1}]$ is a finitely generated projective $R$-module with a flat algebraic connection. 
 Analytification produces a holomorphic vector bundle with flat holomorphic connection over $U(\C)$. Fix $u_0\in U(\C)$.  We have the monodromy representation 
 \[
 \rho : \; \pi_1(U(\C),u_0)  \longrightarrow {\rm  GL } ({\mathcal C}_{\rm ell} ( V)_{u_0} ).
 \]

   We first describe the fundamental group of $ U ( \C )$.
   	For every base point $u_0\in U(\C)$, there is an isomorphism
   	\begin{equation}
   	\pi_1(U(\C),u_0)\cong B_3
   	=\langle\sigma_1,\sigma_2\mid
   	\sigma_1\sigma_2\sigma_1=\sigma_2\sigma_1\sigma_2\rangle.
   	\end{equation}
 We give a sketch of the proof of this fact.
 	Consider the monic cubic
   	\[
   	p_{x,y}(t)=t^3-3xt-2y,\qquad (x,y)=(E_4,E_6).
   	\]
   	Its discriminant equals $108(x^3-y^2)$, so $p_{x,y}$ has three distinct roots precisely when $ ( x , y )\in U$. By Vieta's formulas, $U(\C)$ identifies with the space of unordered triples of distinct complex numbers with zero sum:
   	\[
   	U(\C)\cong F_3^0(\C)/S_3,
   	\quad F_3^0(\C)=\{(z_1,z_2,z_3):z_i\ne z_j,\ \sum_i z_i=0\}.
   	\]
   	Subtracting the center of mass gives an $S_3$-equivariant homeomorphism
   	$F_3(\C)\cong F_3^0(\C)\times\C$, where $S_3$ acts trivially on the second factor. Hence
   	$F_3(\C)/S_3\cong U(\C)\times\C$.
   	The fundamental group of the unordered configuration space $F_3(\C)/S_3$ is, by definition, $B_3$.
   
   \
   
     It is known that  $ \pi_1 (  U , u_0 ) = B_3 $ is a central extension of ${\rm SL}_2 ( \Z )$. We give a proof of this fact.  We 
    consider ${\rm SL}_2 ( \mathbb Z )$-action on $ \mathbb H \times \C^* $ by 
   \[     g \cdot ( \tau , \lambda ) = ( ( a \tau + b ) ( c \tau + b  )^{-1} , ( c \tau + b  )^{-1} \lambda ) , \; \; \; \; \;   g = \left( \begin{matrix} a & b \\ c & d \end{matrix} \right) \]
   This action is a free action.  The orbit space $ ( \mathbb H \times \C^* ) / {\rm SL}_2 ( \Z ) $ is isomorphic to 
    $ U $ by the map
    \[    ( \tau , \lambda ) \mapsto  ( \lambda^4  E_4 ( \tau ) , \lambda^6 E_6 ( \tau  ). \]
    So we have the exact sequence 
 \begin{equation}
     \pi_1 ( \mathbb H \times \C^* ) = \Z \to      \pi_1 ( U ) \to  {\rm SL}_2 ( \Z )      \to 1   .  
 \end{equation}             
   The central element
   \[
   z=(\sigma_1\sigma_2)^3=(\sigma_1\sigma_2\sigma_1)^2
   \]
   generates the center $Z(B_3)= \pi ( \mathbb H \times \C^* ) $.

 \
  
  \begin{theorem}\label{theorem4.4}
   For $M = \mathcal C_{\rm ell} ( V)$, where $V$ is a quasi-lisse VOA,  	the monodromy representation 
    $\rho:B_3\to {\rm GL}(V)$ factors through ${\rm SL}_2(\Z)$.
  \end{theorem}
  
\begin{proof}
  Define the derivation
  \[
  W=4E_4\partial_{E_4}+6E_6\partial_{E_6}.
  \]
  Its algebraic $\C^\times$-action is
  $\lambda\cdot(x,y)=(\lambda^4x,\lambda^6y)$, preserving $U$.
  For $u_0=(a,b)\in U(\C)$, define
  \[
  \gamma_W(t)=(e^{4it}a,e^{6it}b),\qquad 0\le t\le2\pi.
  \]
  Under this action, each root of $p_{x,y}$ scales by $\lambda^2$, since
  \[
  p_{\lambda^4x,\lambda^6y}(\lambda^2z)=\lambda^6p_{x,y}(z).
  \]
  Thus $\gamma_W$ rotates the configuration of roots twice through $2\pi$. One full counterclockwise rotation is the full twist $z$, so
  \[
  [\gamma_W]=z^2.
  \]
  The half-orbit
  \[
  \gamma_{1/2}(t)=(e^{4it}a,e^{6it}b),\qquad 0\le t\le\pi,
  \]
  is also a loop and represents $z$. (Orientations here agree with the standard positive-braid convention.)
  For a homogeneous algebraic section $m\in \mathcal E ( V)_n$, along $\gamma_W$ put
  \[
  s(t)=e^{-int}m|_{\gamma_W(t)}.
  \]
  The velocity is $iW$, so
  \[
  \nabla_{d/dt}s(t)=(-in+in)e^{-int}m|_{\gamma_W(t)}=0.
  \]
  Thus $s$ is horizontal. As $e^{-2\pi in}=1$, parallel transport around $\gamma_W$ fixes the value of $m$ at $u_0$.  Hence $\rho (z^2)=I$.
\end{proof}

  \

 \section { Traces of Vertex Operators and Representations of  $ C_{\rm ell} ( V) $      }\label{section5}

 \subsection { Trace Maps Associated to $V$-modules}\label{section5.1}
 
  Let $V$ be a vertex operator algebra with the central charge of the Virasoro algebra $ c$.
  Let $M $ be a representation of vertex operator algebra $V$. Assume there is a complex number  $h$ and a decomposition 
  $ M = \oplus_{ n = 0 }^\infty M_{ n + h } $ such that ${\rm dim} \, M_{ n + h } < \infty $ and $ L_0 $ acts on $M_{ h + n }$ as
   $( h + n ) {\rm Id} $.  We have an $\C [ E_4 , E_6 ]$-linear map from 
   \[   \C [ E_4 , E_6 ] V \to \C [[ q ]] q^{h -\frac c {24}}  ,   \; \; a \in V \mapsto {\rm Tr}_M (a ) = {\rm Tr}_M  o (a )  q^{ L_0 - \frac c {24}}   \]
   where $ o(a)$ is the component of $Y( a , z )$ that preserves the grading on $M$, for $ a \in V_n $, $ o ( a ) = a_{ ( n -1 )} $.

   It is natural to work on this map using the square-bracket vertex operator algebra structure associated with $V$, see \cite{Zhu}: let 
   \[  Y [ a , z ] = Y ( a , e^{ 2\pi i z } -1  ) e^{ 2\pi i z {\rm deg} (a )} = \sum_{ n = - \infty}^\infty a_{[n]} z^{-n-1} .\]
   $V$ is a vertex operator algebra under the new structure map $ a\mapsto Y [ a , z ]$.  The element $\tilde {\omega} = ( 2\pi i )^2 ( \omega - \frac c { 24}) $
    is the new Virasoso element. We denote this new vertex operator algebra by $ V^{ [ \ ]}$.
    We have 
    \[ a_{[m]} = ( 2\pi i )^{-m-1} Res_z ( Y(a , z ) ( \ln ( 1 + z ))^m  (1 + z)^{ deg (a )-1 } ).  \]
     The results in \cite{Zhu} imply that this maps factors through ${\mathcal C}_{\rm ell} ( V^{ [ \ ]}) $.

 \begin{prop}\label{prop5.1}
 (1).	The map ${\rm Tr}_M$ factors through $C_{\rm ell} ( V^{ [ \ ]}) $. 
 	\newline (2).  For $ a\in V_k$, 
 	\begin{equation}
 	{\rm Tr}_M ( D a ) = \left(  q \frac d { d q } - \frac { k } { 1 2 } E_2 (\tau ) \right)  	{\rm Tr}_M ( a ) . 
 	\end{equation}   
 	\newline (3).   If $V$ is quasi-lisse,  $ {\rm Tr}_M ( a ) $ converges to an analytic function on $ 0 < |q | < 1 $, therefore
 	 the map ${\rm Tr}_M$ is $\C [ E_4 , E_6 ]$-linear map from 
 	$  \C [ E_4 , E_6 ] V $ to  $ {\rm Hol}  ( \mathbb H  )$, the space of holomorphic functions on the upper half space $\mathbb H $.    
 \end{prop}
 
   \begin{proof}
   	 (1) follows from the Proposition 4.3.6 \cite{Zhu}.  (2) follows from Proposition 4.3.6 in\cite{Zhu}.
   	(3) uses (2), the detailed proof can be found in \cite{AK}.
   \end{proof}
   
   \

  \subsection {Coinvariants for Nodal and Cuspidal Cubic Curves  }\label{s5.2}

  For a point $ p = ( a , b ) \in \C^2 $, we have the evaluation homomorphism 
  \[  ev_p : \;  \C [ E_4 , E_6 ] \to   \C \qquad    E_4 \mapsto    a , \,  E_6 \mapsto    b  \]
  The fiber 
  \[   C_{\rm ell} ( V)_p = C_{\rm ell} ( V) \otimes_{ \C [ E_4 , E_6 ] ,  ev_p } \C \]
  is the space of the coinvariant for the curve $ \pi^{-1} ( p ) $ with $V$ assigned at the marked point $O $. 
  Let $A(V)$ denote Zhu's algebra of $V$ and $ R_V$ denote Zhu's Poisson algebra of $V$.  The coinvariant on cubic nodal and cuspidal curves are related 
  to $A(V)$ and $R_V$ as follows.

   \begin{prop}\label{prop5.2}
   For  $ p = (a , b)\in \C^2 $,	 if  $ \pi^{-1} p $ is a  nodal cubic curve, that is  $a^3 - b ^2 = 0 $, $ ( a , b ) \ne ( 0 , 0 )$, then the fiber 
   	 $ C_{\rm ell} ( V^{ [ \ ]})_p $ is isomorphic to $ A(V) / [ A( V ) , A(V) ] $.  
   	 If  $ \pi^{-1} p $ is a cuspidal cubic curve, that is , $ ( a , b )= ( 0 , 0 )$, then the fiber  $ C_{\rm ell} ( V^{ [ \ ]})_p $ is isomorphic to
    to $ R_V /[\{ R_V  , R_V \}  $. 
   \end{prop}

 The space of coinvariants can be identified with the zeroth chiral homology. In the nodal degeneration of an elliptic curve, the zeroth chiral homology is identified with the zeroth Hochschild homology of  $A(V)$, while the first chiral homology is closely related to the first Hochschild homology of $A(V)$ \cite{EH}. Our proposition suggests investigating the first chiral homology of a cuspidal cubic curve and its relationship with the first Poisson homology of the $C_2$-algebra $R_V$.
 
  For a rational vertex operator algebra with $C_2$-condition, we have 
  
\begin{theorem}\label{theorem5.3}
 If $V$ is a rational vertex operator algebra with $C_2$-condition, then $ C_{\rm ell} (  V  )_{\rm tor}  $ is annihilated by $ E_4^m $ and $ E_6^m$ for some $m$.
 And $ {\rm dim} \, C_{\rm ell} (  V  )_{\rm tor} < \infty $. 
 \end{theorem}

\begin{proof}
 $ C_{\rm ell} (  V  )_{\rm tor}  $ is a finite over $\C [ E_4 , E_6 ]$.
 Consider the sheaf ${\rm dim} \, C_{\rm ell} (  V  )_{\rm tor}^{\rm sim} $ on ${\rm Spec} \C [ E_4 , E_6 ]$, our conditions imply that it is supported at $( 0 , 0 )$.
 This proves the results. 
\end{proof}

  \
  
   \subsection { $\C [ E_4 , E_6 ] \la D, W \ra $-Module Formed by Holomorphic Functions on the Upper Half Plane $\mathbb H$  }\label{s5.3}
  
  To interpret the trace map $ {\rm Tr}_M $ as a homomorphism of  $\C [ E_4 , E_6 ] \la D, W \ra $-Modules, we need to 
  regard the images $	{\rm Tr}_M $ of different homogeneous components as distinct spaces.    
  Accordingly, for each $ n \geq 0 $, we  we introduce a copy $ {\rm Hol}  ( \mathbb H  )_n $ of $ {\rm Hol}  ( \mathbb H  )$,
  and set 
  \[  \mathcal H  =  \oplus_{ n = 0 }^\infty  {\rm Hol}  ( \mathbb H  )_n .\]
 We define a $\C [ E_4 , E_6 ] \la D, W \ra $-module on $\mathcal H$ as follows.
  For $ f (\tau ) \in C [ E_4 , E_6 ] $ with weight $ k $,  $ \phi ( \tau ) \in  {\rm Hol}  ( \mathbb H  )_n $, 
  the scalar multiplication $ f \cdot \phi $ is the function  $ f ( \tau ) \phi ( \tau ) $ considered as an element in $ {\rm Hol}  ( \mathbb H  )_{k + n }$.
  This gives  $ \mathcal H $ a structure of graded $ \C [ E_4 , E_6 ]$-module. 
  $ W$ acts on  ${\rm Hol}  ( \mathbb H  )_n $ as $ n \, Id $.  The action of $ D $ on $ {\rm Hol}  ( \mathbb H  )_n $
   is the map
   \[    D : \, {\rm Hol}  ( \mathbb H  )_n \to {\rm Hol}  ( \mathbb H  )_{n+2 } \]
   given by 
   \begin{equation}
      D =  ( 2\pi i )^{-1} \frac d { d\tau } - \frac { n } { 1 2 } E_2 (\tau ).
   \end{equation} 
  It is easy to check that $ \mathcal H $ is a $\C [ E_4 , E_6 ] \la D, W \ra $-module under the above actions.
  The following result is a direct corollary of Proposition \ref{prop5.1}.

  \

  \begin{prop}\label{prop5.4}
  We define the grading preserving liner
  map 
  \[  {\rm T}_M :   \, C_{\rm ell} ( V ) = \oplus_{ n = 0 } {C}_{\rm ell} ( V^{ [\ ]} )_n   \to \mathcal H = \oplus_{ n= 0}^\infty{\rm Hol}  ( \mathbb H  )_n \]
  by
  \[  {\rm T}_M ( a f ( \tau )) = f (\tau) \, {\rm Tr}_M o ( a ) q^{ L_0 - \frac c { 24} } .\]
  Then $ {\rm T}_M$ is a homomorphism of  $\C [ E_4 , E_6 ] \la D, W \ra $-modules.
  \end{prop}

  \
  
  \subsection { ${\rm SL}_2 ( \Z )$-action   }\label{s5.5}

  \begin{prop}\label{prop5.5}
  	Let $SL_2 ( \Z )$ acts on  $ \mathcal H $ as follows. 
  	for 
  	\[ g = \left( \begin{matrix}  a & b \\ c & d \end{matrix} \right) \in SL_2 ( \Z ), \; \;  f (\tau ) \in  {\rm Hol } ( \mathbb H )_n \]
  	\[  ( g \cdot f ) (\tau )  = ( c \tau + d )^{-n }	f \left( \frac { a \tau + b } { c\tau + d }   \right) \in  {\rm Hol } ( \mathbb H )_n\]
  	This action is a right action and it preserves the $\C [ E_4 , E_6 ] \la D, W \ra $-structure.
  \end{prop}

  \
  
\begin{proof}   If $ f \in \C [ E_4 , E_6] $ is a module form of weight $ k $,   $  \phi (\tau ) \in {\rm Hol} ( \mathbb H  )_n $, then $ f \phi \in 
		{\rm Hol} ( \mathbb H)_{k+n} $,
  \[  g \cdot ( f  \phi ) (\tau )  =   ( c \tau + d )^{ - k -n } f \left( \frac { a \tau + b } { c\tau + d }   \right) \phi \left( \frac { a \tau + b } { c\tau + d }   \right)
  = f (\tau )   ( c \tau + d )^{  -n } \phi \left( \frac { a \tau + b } { c\tau + d }   \right) = f ( g \cdot \phi ) .\]
  \begin{eqnarray*}
  	g \cdot ( D \phi ) (\tau ) &=&   ( c \tau + d )^{  -n -2 }   ( D \phi )\left( \frac { a \tau + b } { c\tau + d }   \right)\\
  	&=&   ( c \tau + d )^{  -n -2 }  \left( (2 \pi i)^{-1} (\partial_\tau \phi ) \left( \frac { a \tau + b } { c\tau + d }   \right) 
  	- \frac n {12}     E_2 \left( \frac { a \tau + b } { c\tau + d }   \right)\phi  \left( \frac { a \tau + b } { c\tau + d }   \right) \right)  
  \end{eqnarray*}

  \begin{eqnarray*}
  	D   ( g \cdot \phi ) ( \tau )  &=&    (  (2\pi i)^{-1} \partial_\tau -\frac n {12} E_2 ( \tau )  ) \left(   ( c \tau + d )^{-n}\phi  \left( \frac { a \tau + b } { c\tau + d }   \right)
  	\right) \\
  	&=&  (2\pi i)^{-1} (c\tau + d )^{-n-2}  (\partial_\tau \phi ) \left( \frac { a \tau + b } { c\tau + d }   \right) - (2\pi i)^{-1} n c ( c \tau + d )^{-n-1}\phi  \left( \frac { a \tau + b } { c\tau + d }   \right) \\
  	&& \; \; \; \; \;  -\frac {n} {12} E_2 ( \tau )    ( c \tau + d )^{-n}\phi  \left( \frac { a \tau + b } { c\tau + d }   \right)
  \end{eqnarray*}
  We see that 
  \[     g \cdot ( D \phi ) = D   ( g \cdot\phi  ) \]
  following the identity
  \[   E_2 \left( \frac {a \tau + b } { c\tau + d }  \right) = ( c\tau + d )^2 E_2 (\tau ) -6 i \pi^{-1} c ( c\tau + d ) \]
  which is equivalent to (\ref{2.6})  
 \end{proof}
  
  \

 As a consequence of Proposition \ref{prop5.5},  $ SL_2 ( \Z )$  acts from the right on the space of homomorphisms by post-composition
\begin{equation}
  {\rm Hom}_{ \C [ E_4 , E_6 ] [D , W  ]  } (    { C}_{\rm ell} ( V^{ [ \ ] }) ,  \, \mathcal H  ).
  \end{equation}
We call this representation ${\rm SL}_2 (\Z)$-representation on elliptic homomorphisms of 
 $V$.
  If $V$ is rational and satisfies the $C_2$-condition,  the results in \cite{Zhu} implies that 
  $ T_M $, as $M$ runs through the  irreducible $V$-modules, gives a basis for this space.

  \
  
 \begin{theorem}\label{theorem5.6}
 If $V$ is a rational vertex operator algebra and satisfies the $C_2$-condition, then 	the monodromy representation of ${\rm SL_2} ( \Z )$ and  ${\rm SL}_2 (\Z)$-representation on elliptic homomorphisms
 	are dual each other. 
  \end{theorem}
  
  If $V$ is quasi-lisse, we can prove there are $SL_2 ( \Z )$-invariant pairing of the two representations.

\


\begin{thebibliography}{aa}

\bibitem[AN]{AN} T. Abe and K.Nagatomo. Finiteness of conformal blocks over compact Riemann
surfaces. Osaka J. Math., 40(2):375–391, 2003.

\bibitem[AK]{AK} Tomoyuki Arakawa, Kazuya Kawasetsu, 
\emph{Quasi-lisse Vertex Algebras and Modular Linear Differential Equations}
Published In: Lie Groups, Geometry, and Representation Theory (A Tribute to the Life and Work of Bertram Kostant), Progress in Mathematics, Vol. 326 
Pages: 41--57, 2018 s

 \bibitem[BD]{BD}
 A.~Beilinson and V.~Drinfeld,
 \emph{Chiral Algebras},
 American Mathematical Society Colloquium Publications, vol.~51,
 American Mathematical Society, Providence, RI, 2004.
 \url{https://doi.org/10.1090/coll/051}.
 
  \bibitem[EH]{EH}
 J. van Ekeren and R. Heluani, 
 \emph{Chiral Homology of Elliptic Curves and the Zhu Algebra},  Communications in Mathematical Physics, 386 (2021), no. 1, pp. 495–550
 

 
 
 \bibitem[DGT21]{DGT21}
 C.~Damiolini, A.~Gibney, and N.~Tarasca,
 ``Conformal blocks from vertex algebras and their connections on $\overline{\mathcal M}_{g,n}$,''
 \emph{Geometry \& Topology} \textbf{25} (2021), no.~5, 2235--2286.

 
 \bibitem[DGT24]{DGT24}
 C.~Damiolini, A.~Gibney, and N.~Tarasca,
 ``On factorization and vector bundles of conformal blocks from vertex algebras,''
 \emph{Annales scientifiques de l'\'Ecole normale sup\'erieure} (4)
 \textbf{57} (2024), no.~1, 241--292.

 

\bibitem[FBZ]{FBZ}
E.~Frenkel and D.~Ben-Zvi,
\emph{Vertex Algebras and Algebraic Curves},
2nd ed., Mathematical Surveys and Monographs, vol.~88,
American Mathematical Society, Providence, RI, 2004.


 
\bibitem[NT]{NT}  K. Nagatomo and A. Tsuchiya. \emph{Conformal field theories associated to regular
 chiral vertex operator algebras. I}. Theories over the projective line. Duke Math.
 J.,
 



\bibitem[Sil]{Sil} J. H. Silverman,  Advanced Topics in the Arithmetic of Elliptic Curves,  GTM 151, Springer-Verlag, 1994

 
 \bibitem[TUY]{TUY}
 A.~Tsuchiya, K.~Ueno, and Y.~Yamada,
 ``Conformal field theory on universal family of stable curves with gauge symmetries,''
 \emph{Advanced Studies in Pure Mathematics} \textbf{19} (1989), 459--566.
 \url{https://doi.org/10.2969/aspm/01910459}.
 

 
 \bibitem[Z]{Zhu}
 Y.~Zhu,
 ``Modular invariance of characters of vertex operator algebras,''
 \emph{Journal of the American Mathematical Society} \textbf{9} (1996), no.~1, 237--302.
 \url{https://doi.org/10.1090/S0894-0347-96-00182-8}.
 

\end{thebibliography}
\end{document}